\documentclass[11pt,a4paper]{scrartcl}

\makeatletter
\DeclareOldFontCommand{\rm}{\normalfont\rmfamily}{\mathrm}
\DeclareOldFontCommand{\sf}{\normalfont\sffamily}{\mathsf}
\DeclareOldFontCommand{\tt}{\normalfont\ttfamily}{\mathtt}
\DeclareOldFontCommand{\bf}{\normalfont\bfseries}{\mathbf}
\DeclareOldFontCommand{\it}{\normalfont\itshape}{\mathit}
\DeclareOldFontCommand{\sl}{\normalfont\slshape}{\@nomath\sl}
\DeclareOldFontCommand{\sc}{\normalfont\scshape}{\@nomath\sc}
\makeatother

\usepackage{a4wide}
\usepackage{latexsym}
\usepackage[USenglish]{babel}
\usepackage{appendix}
\usepackage[utf8]{inputenc}
\usepackage[T1]{fontenc}
\usepackage{amssymb}
\usepackage{hyperref}
\usepackage[amsthm,thmmarks]{ntheorem}
\usepackage{fixmath}

\usepackage{mathrsfs}
\usepackage{graphicx}
\usepackage{amsmath}
\usepackage{mathdots}
\usepackage{amsfonts}
\usepackage{amsxtra}
\usepackage{url}
\usepackage{mathptmx}
\usepackage{helvet}
\usepackage{stmaryrd}
\usepackage{dsfont}
\usepackage{setspace}
\usepackage{enumerate}
\usepackage{verbatim}
\usepackage{rotating}
\usepackage{color}
\usepackage{cleveref}
\usepackage{cancel}
\usepackage{epstopdf}
\usepackage[tablewithin=section]{caption}
\usepackage{MathDef}

\def\dif{\mathrm{d}}

\def\diag{\text{diag}}
\DeclareMathOperator{\res}{Res}

\theoremstyle{plain}
\newtheorem{theorem}{Theorem}[section]
\newtheorem{lemma}[theorem]{Lemma}
\newtheorem{proposition}[theorem]{Proposition}
\newtheorem{definition}[theorem]{Definition}
\newtheorem{corollary}[theorem]{Corollary}
\newtheorem{assumptionletter}{{\textbf{Assumption}}}

\theoremstyle{definition}
\newtheorem{example}[theorem]{Example}
\newtheorem{remark}[theorem]{Remark}

\numberwithin{equation}{section}
\allowdisplaybreaks

\newcommand{\VF}[1]{{\color{blue} #1}}

\title{Factorization and discrete-time representation  of multivariate \\ CARMA processes}

\author{Vicky Fasen-Hartmann \setcounter{footnote}{1}\thanks{Institute of Stochastics, Englerstra{\ss}e 2,
D-76131 Karlsruhe, Germany, email: vicky.fasen@kit.edu} \label{fnote}%\emph{Email:}
 \and Markus Scholz \thanks{Allianz Lebensversichung-AG, Reinsburgstra{\ss}e 19, D-70197 Stuttgart, Germany. }
 \thanks{Financial
support by the Deutsche Forschungsgemeinschaft through the research
grant FA 809/2-2 is gratefully acknowledged.}
}

\date{}

\begin{document}

%%%%%%%%%%%%%%%%%  Titel und Autoreninformationen  %%%%%%%%%%%%%%%%%%%%%%%%%%%%%%%%
\maketitle
%
%%%%%%%%%%%%%%%%%  Abstract und AMS Classification  %%%%%%%%%%%%%%%%%%%%%%%%%%%%%%%
\begin{abstract} \vspace*{-1cm}
In this paper  we show that stationary and non-stationary
multivariate continuous-time ARMA (MCARMA) processes have the representation as
a sum of  multivariate complex-valued Ornstein-Uhlenbeck processes under some
mild assumptions. The proof benefits from properties of
rational matrix polynomials. A conclusion is an alternative description of
the autocovariance function of a stationary MCARMA process.
Moreover, that representation is used
to show that the discrete-time sampled MCARMA$(p,q)$  process is a weak VARMA$(p,p-1)$
process if second
moments exist. That result complements the weak VARMA$(p,p-1)$ representation
derived in Chambers and Thornton~\cite{Chambers:Thornton:2012}. In particular,
it relates the right solvents of the autoregressive polynomial of the MCARMA process
to the right solvents of the autoregressive polynomial of the VARMA process;
in the one-dimensional case the right solvents are the zeros of the autoregressive polynomial.
Finally, a factorization of the sample autocovariance function of the noise sequence  is presented
 which is useful for
statistical inference. %Analog results for integrated MCARMA processes are derived as well.
\end{abstract}

\noindent
\begin{tabbing}
\\ \> Secondary:  	62M86,  	60G10     	
\end{tabbing}

\vspace{0.2cm}\noindent\emph{Keywords:} autocovariance function, latent root,  matrix polynomial, MCARMA process, Ornstein-Uhlenbeck process, rational matrix function, right solvent, VARMA process.

%%%%%%%%%%%%%%%%%  Einleitung  %%%%%%%%%%%%%%%%%%%%%%%%%%%%%%%%%%%%%%%%%%%%%%%%%
\section{Introduction}
%
%
%%%%%%%%%%%%%%%%%  Literaturverzeichnis  %%%%%%%%%%%%%%%%%%%%%%%%%%%

A multivariate continuous-time ARMA (MCARMA) process is a continuous-time version of the well-known
vector ARMA (VARMA) process in discrete time. They are applied in diversified fields as, e.g., signal processing and control (cf. \cite{GarnierWang2008,LarssonMossbergSoederstroem2006}),
high-frequency financial econometrics (cf. \cite{Todorov2009})  and financial mathematics (cf. \cite{BenthKoekebakkerZakamouline2010}).
The driving process of a MCARMA process is a L\'evy process $L=(L(t))_{t\geq 0}$ which is
an $\R^m$-valued stochastic process with $L(0)=0_m$ $\mathbb{P}$-a.s., stationary and independent increments and c\`adl\`ag sample paths. The idea is then that a $d$-dimensional MCARMA$(p,q)$ process  ($p>q$ positive integers)
 is  the solution of the stochastic differential equation
\begin{eqnarray} \label{eq1.1}
     {A}(D)Y(t)={B}(D)D L(t) \quad \mbox{ for } t\geq 0,
\end{eqnarray}
where $D$ is the differential operator with respect to $t$,
\begin{equation} \label{Pol}
    {A}(\lambda):= I_{d}\lambda^p+A_1\lambda^{p-1}+\ldots+A_p \quad \text{ and }\quad
    {B}(\lambda):=B_0\lambda^q+\ldots+ B_{q-1}\lambda+B_q
\end{equation}
 is the autoregressive and the moving average polynomial, respectively with $A_1,\ldots, A_p\in \R^{d\times d}$  and
$B_0,\ldots,B_q\in \R^{d\times m}$.
The matrix %$0_{d\times s}$
%is the zero matrix in $\R^{d\times s}$ and
$I_{d}$ denotes the
$d\times d$-dimensional identity matrix and $0_{d\times m}$ denotes a $d\times m$-dimensional matrix whose entries are all zero in the following. In contrast, in discrete time the differential operator is replaced by the backshift operator
and the differential of the Lévy process $D L(t)$ by a weak white noise.
 Since a Lévy process is not differentiable, the question arises what is the formal
definition of a MCARMA process. We can interpret \eqref{eq1.1} via linear continuous-time state space models as in Marquardt and Stelzer \cite{MarquardtStelzer2007}.
Therefore,
define
\beam \label{def:A:star}
    A^*:=\left(\begin{array}{ccccc}
        0_{d\times d} & I_{d} & 0_{d\times d} & \cdots & 0_{d\times d}\\
        0_{d\times d} & 0_{d\times d} & I_{d} & \ddots & \vdots \\
        \vdots & & \ddots & \ddots & 0_{d\times d}\\
        0_{d\times d} & \cdots & \cdots & 0_{d\times d} & I_{d}\\
        -A_p & -A_{p-1} & \cdots & \cdots & -A_1
    \end{array}\right) \in \R^{pd\times pd},
\eeam
$C^*:=(I_{d},0_{d\times d},\ldots,0_{d\times d}) \in \R^{d\times pd}$ and $B^*:=(\beta_1^\mathsf{T} \cdots \beta_p^\mathsf{T})^\mathsf{T}\in \R^{pd\times m}$ with $\beta_1:=\ldots:=\beta_{p-q-1}:=0_{d\times m}$ and
\beao
    \beta_{p-j}:=-\sum_{i=1}^{p-j-1}A_i \beta_{p-j-i}+B_{q-j},  \quad j=0,\ldots,q.
\eeao
Then the $\R^{d}$-valued \textsl{ MCARMA$(p,q)$ process} $Y:=(Y(t))_{t\geq 0}$ is  defined by the state space
equation
\begin{align}   \label{state:space}
\begin{split}
    Y(t)&=C^* X(t) \quad \text{ and } \quad %\quad \mbox{ for } t\geq 0,\\
\dd X(t)=A^* X(t)\,\dd t+B^* \,\dd L(t).
\end{split}
\end{align}
Interesting is that if we define
\beam \label{A hash}
    A^{\#}=\left(\begin{array}{ccccc}
        I_{d} & 0_{d\times d}   & \cdots & 0_{d\times d}\\
        A_1 & I_{d} & \ddots &  \vdots \\
       % A_2 &  \ddots &  \ddots & \ddots  &  \vdots \\
         \vdots &   \ddots  & \ddots & 0_{d\times d}\\
        A_{p-1} &  \cdots & A_1 & I_{d}
    \end{array}\right) \in \R^{pd\times pd}, \quad
    B^{\#}=\left(\begin{array}{c} 0_{(p-(q+1))d\times m}\\
    B_0\\
    \vdots\\
    B_{q}\\
    \end{array}\right)\in \R^{pd\times m},
\eeam
then
\beam \label{2}
    A^{\#}B^{*} =B^{\#}.
\eeam
%This is a very basic relation for this paper because a consequence of  Kailath~\cite{Kailath}
%(cf. Levya-Ramos \cite[Section 6, eq. (6.3)]{Levya-Ramos1991})
%is that
%\beam \label{2.2}
%    C^*(\lambda I_{pd} -A^*)^{-1}B^*=A(\lambda)^{-1}B(\lambda).
%\eeam
%Sato and Yamazato~\cite[Theorem 4.1 and Theorem 4.2]{Sato:Yamazato:1984} derived necessary and sufficient conditions
%for the existence of the multivariate Ornstein-Uhlenbeck process $X$. A direct consequence from this
%is that if the spectrum  $\sigma(A^*):=\{\lambda\in\C:\det(\lambda I_{pd}-A^*)=0\}\subseteq (-\infty,0)+i\R$ and %$\E[\log(\max(1,\|L(1)\|))]<\infty$ then
%\beam \label{MCARMA}
%    Y(t)=\int_{- \infty}^{\infty}g(t-s)L(\dif s), \quad t\in\R,
%\eeam
%with $g(s)= C^*\mathrm{e}^{A^*s}B^* \mathds{1}_{\left[0,\infty\right)}(s)$ for $s\in\R$
%is a stationary MCARMA process and solves the state space equation \eqref{state:space}.
The class of MCARMA processes is very rich.
Under the constrain of finite second moments,
Schlemm and Stelzer \cite[Corollary 3.4]{SchlemmStelzer2012}  show that the class of stationary MCARMA processes
and the class of stationary state space models are equivalent (see Fasen-Hartmann and Scholz~\cite{FasenScholz1} for cointegrated MCARMA processes).

The aim of the paper is to present sufficient criteria for  stationary and non-stationary
MCARMA$(p,q)$ processes
to have a representation as  a sum of $p$ multivariate Ornstein-Uhlenbeck processes (which are
 MCAR(1) $=$  MCARMA(1,0) processes).
For $d=1$, it is well known that if the zeros $r_1,\ldots,r_p$ of $A(\lambda)$ are distinct
    and have a strictly negative real parts that
 then
 \beam \label{SOU}
    Y(t)=\sum_{k=1}^pY_k(t) \quad \text{ with } \quad
    Y_k(t)=\int_{-\infty}^t\e^{r_k (t-u)}\,\frac{B(r_k)}{A'(r_k)}\,\dif L(u)
 \eeam
 is a stationary solution of the state space model \eqref{state:space} and hence, a CARMA process (see Brockwell,
 Davis and Yang~\cite[Proposition 2]{BrockwellDavisYang2011}). Note that $B(r_k)/A'(r_k)$ is the residue of $A(\lambda)^{-1}B(\lambda)$ at $r_k$.
 In the present paper we extend this result to the multivariate setup for both stationary and non-stationary
 MCARMA processes. The zero $r_k$ of $A(\lambda)$ in the one-dimensional case
 is replaced by a $d\times d$ matrix $R_k$, which is as well
 a kind of multivariate ''zeros'' of $A(\lambda)$, the so called right solvent satisfying
 $A_R(R_k):=I_dR_k^p+A_1R_k^{p-1}+\ldots+A_p=0_{d\times d}$. The result is derived in \Cref{thmBlockDecompStateSpaceMod}.
 Essential for our proof are basic principles from rational matrix polynomials coming from linear algebra
 which are not necessary in dimension $d=1$. A main feature is that we have a sum of multivariate
 Ornstein-Uhlenbeck processes and not only some linear combinations of multivariate
 Ornstein-Uhlenbeck processes. Since
   matrix multiplication is not commutative this is not trivial. That is different to the one-dimensional case where
 any linear combination of stationary Ornstein-Uhlenbeck processes is as
 well a sum of stationary Ornstein-Uhlenbeck processes.
A straightforward consequence of our result is an alternative representation of the autocovariance function
of a stationary MCARMA process in \Cref{Propo_acf}.

Although we consider in this paper a continuous-time model, the corresponding discrete-time models are of special interest. The reason for this is that despite having a continuous-time model, one often observes  the process only at discrete time points
as, e.g, in the context of high-frequency data.
Hence, we use the representation of a MCARMA$(p,q)$ process as a sum of $p$ multivariate Ornstein-Uhlenbeck
processes to derive a vector-valued ARMA (VARMA)$(p,p-1)$ representation
for the low frequency sampled MCARMA process $(Y(nh))_{n\in\N}$ ($h>0$ fixed) in \Cref{thmBlockweakVARMArep}.
For the proof of this theorem a representation of $Y$ as a  linear combination of multivariate Ornstein-Uhlenbeck processes
is not sufficient. The statement is a direct extension of the
ARMA$(p,p-1)$ representation of discretely sampled CARMA processes in
Brockwell, Davis and Yang \cite[Proposition 3]{BrockwellDavisYang2011} whose
 autoregressive polynomial $\prod_{k=1}^p(1-\lambda \e^{r_k h})$ of the ARMA representation has as zeros $\e^{-r_1 h},\ldots,\e^{-r_p h}$.
 In analogy, in the multivariate setup of this paper
 the autoregressive polynomial of the VARMA representation has right solvents  $\e^{-R_1 h},\ldots,\e^{-R_p h}$.

 In the econometric literature,  the VARMA$(p,p-1)$ representation of a discretely sampled \linebreak
MCARMA process is well-known, see, e.g., Chambers and Thornton~\cite[Corollary 1]{Chambers:Thornton:2012};
a nice overview on this topic is presented in Chambers, McCrorie and Thornton \cite{ChambersMcCrorieThornton}.
In contrast to us, Chambers and Thornton~\cite{Chambers:Thornton:2012} assume some
kind of observability  and controllability conditions
on submatrices of $\e^{A^{**}}$, where $A^{**}$ is constructed form $A^*$
by reflecting the entries of $A^*$ at the diagonal from the left lower corner to the right upper
corner. There, the coefficients of the autoregressive polynomial in the VARMA representation are complicated functions of these submatrices.
  The current paper presents an alternative and simpler representation of the VARMA parameters and in particular, % statement gives  more inside into the structure of the autoregressive polynomial in the VARMA representation
 it connects the autoregressive polynomial in the MCARMA representation
 to the autoregressive polynomial in the VARMA representation due to the solvents.
 Our proof is an alternative proof requiring only assumptions
 on the right solvents of $A(\lambda)$.
%In particular, we think that it might be easier to identify parameters from the MCARMA process VARMA representation
%In the case of CARMA processes details can be found in Brockwell and Lindner~\cite{Brockwell:Lindner:2019}
%and  Brockwell, Davis and Yang \cite[Proposition 2]{BrockwellDavisYang2011}.
% Brockwell, Davis and Yang \cite[Proposition 2]{BrockwellDavisYang2011}
%already showed that if the eigenvalues of $A^*$ are distinct then a
%CARMA process has a  decomposition  into a sum of $p$ dependent Ornstein-Uhlenbeck  processes.
%Based on this they first developed a weak VARMA$(p,p-1)$ representation of the
%discretely-sampled CARMA process and then an estimation method for the CARMA parameters.
%Our main statements, \Cref{thmBlockDecompStateSpaceMod} and \Cref{thmBlockweakVARMArep}, are multivariate
%extensions of their results.
In the multivariate setting, Schlemm and Stelzer~\cite[Proposition 5.1]{SchlemmStelzer2012}
proved that a MCARMA process has a representation
 as a multivariate linear combination of $pd$ dependent one-dimensional Ornstein-Uhlenbeck processes.
In the present paper, we will have multivariate Ornstein-Uhlenbeck processes and $p$ instead of $pd$
Ornstein-Uhlenbeck processes.  %Schlemm and Stelzer~\cite{SchlemmStelzer2012}
%apply their proposition  to derive a weak VARMA($pd$,$pd-1$) representation of the discrete-time
%sampled MCARMA process and they analyze the mixing properties of the weak white noise.

Similarly, as in the above mentioned  papers our conclusions are advantageously  for statistical
inference of MCARMA processes.
Brockwell and Lindner~\cite{Brockwell:Lindner:2019} use the representation
\eqref{SOU} to solve both the sampling and the embedding problem for CARMA processes. In the first
case, they deduce the explicit parameters of the ARMA representation of $(Y(nh))_{n\in\N}$.
In the second case, they present conditions under which an ARMA$(p,q)$ process
can be embedded in a CARMA$(p,p-1)$ process. Therefore, we think
that our results might be helpful for a multivariate version of the sampling and embedding problem as well.
But this is outside the scope of the present paper.
Moreover, our findings are helpful to derive
probabilistic properties of a MCARMA process.
Brockwell and Lindner~\cite{Brockwell2009}, for example, use the ARMA$(p,p-1)$ representation of a
CARMA process to derive necessary and sufficient conditions for the existence of a CARMA process.

The paper is structured on the following way. In \Cref{preliminaries} we present preliminary results on matrix polynomials
and rational matrix polynomials which lay the background for the upcoming results. The main results of the paper are given
in \Cref{ch:chapter3:sec:DecompMCARMA}. %Finally, we end with a conclusion in

\section{Preliminaries} \label{preliminaries}

In this section, we review main results on matrix polynomials and rational matrix functions. References about matrix analysis and matrix polynomials are, e.g., the textbooks of Bernstein~\cite{Bernstein2009}, Horn and Johnson \cite{HornJohnson2013} and  Kailath~\cite{Kailath}.% and Gohberg et al. \cite{Gohberg1982}.
The aim is to receive matrix valued ''roots''
of a matrix polynomial  which help to define linear factors of a matrix polynomial. % as in the one-dimensional case.
 However, a challenge is that
 there does not exist
the Fundamental Theorem of Algebra for matrix polynomials and matrix multiplication is not commutative.
%Since in general matrix multiplication is not commutative, it is a difference if we multiply a matrix form the left or the right. In this paper, we consider only the case where the matrix is multiplied from the left. In the other case the definitions and results are analogous.

\begin{definition}
\label{defLambdaMatrixMatrixPolynomial}
$\mbox{}$
\begin{enumerate}[(a)]
\item
A \textsl{$\lambda$-matrix $A:\C\to \C^{d\times m}$ of degree $p$ and order $(d,m)$}  is defined as
\begin{align*}
&A(\lambda)=A_0\lambda^p+A_1\lambda^{p-1}+\ldots+A_{p-1}\lambda+A_p, \quad \lambda \in\C,
\end{align*}
where $A_k\in \C^{d\times m}$ for $k=0,\ldots,p$. If additionally, $d=m$ we say shortly
that $A(\lambda)$ is of degree $p$ and order $d$, and define the spectrum of $A(\lambda)$
as $\sigma(A(\cdot)):=\{\lambda\in\C:\det(A(\lambda))=0\}$.  If $\sigma(A(\cdot))$ lies in the complement of the
closed unit disc, then $A(\lambda)$ is called
\textsl{Schur-stable}.
The $\lambda$-matrix $A(\lambda)$
is called  \textsl{ monic $\lambda$-matrix of degree $p$ and order $d$}
if  $A_0:=I_d$.
\item
Let $Z\in\C^{d\times d}$ and $d=m$. Then the \textsl{right matrix polynomial} $A_R: \C^{d\times d}\to \C^{d\times d}$
of the $\lambda$-matrix $A(\lambda)$ is defined as
\begin{align*}
A_R(Z):=A_0Z^p+A_1Z^{p-1}+\ldots+A_{p-1}Z+A_p.
\end{align*}
%\item[(c)] An \textsl{$q^{th}$-degree, $(d,m)^{th}$-order  $\lambda$-matrix} $B:\C\to \C^{d\times m}$ is given by
%\begin{align*}
%&B(\lambda)=B_0\lambda^q+B_1\lambda^{q-1}+\ldots+B_{q-1}\lambda+B_q, \quad \lambda \in\C,
%\end{align*}
%where $B_k\in \C^{d\times m}$ for $k=0,\ldots,q$.
\end{enumerate}
%and analogously the \textsl{left matrix polynomial} $A_L$ is given by
%$$A_L(Z)=Z^nA_0+Z^{n-1}A_1+\ldots+ZA_{n-1}+A_n.$$
\end{definition}
%In this paper we restrict to right matrix polynomials but similarly we could also define the left matrix
%polynomial $A_L(Z)=Z^pA_0+Z^{p-1}A_1+\ldots+ZA_{p-1}+A_p.$ Since matrix multiplication is not commutative the left
%and the right matrix polynomial differ.
%%-----------------------------------------------------------
%By considering the determinant of a $\lambda$-matrix, we obtain a one-dimensional polynomial with ''classical'' roots. The roots of this determinant  and the latent vectors related to these roots, which we define in the following, are important for the multivariate extension of roots of univariate  polynomials.
%%-----------------------------------------------------------
Next, we extend the definition of a root to the matrix polynomial case.
\begin{definition}
\label{defMultipleRoot}
For a  monic $\lambda$-matrix $A(\lambda)$ of degree $p$ and order $d$ we define
\begin{equation*}
A^{(k)}(\lambda):=\frac{d^k}{d\lambda^k}A(\lambda), \quad k=1,\ldots p.
\end{equation*}
A matrix $R\in\C^{d\times d}$ is defined to be a \textsl{right solvent of $A(\lambda)$ with multiplicity $\nu\in\{1,\ldots,p\}$} if
\begin{align*}
%\label{eqCommonRootEquations}
A_R(R)=0_{d\times d},\quad A^{(1)}_R(R)=0_{d\times d},\quad\ldots,\quad A^{(\nu-1)}_R(R)=0_{d\times d}\quad \text{and} \quad A^{(\nu)}_R(R)\not=0_{d\times d}.
\end{align*}
%Thus,  $A(\lambda)$ has the representation as
%$
%\label{EqFactorizationMultipleRoot}
%A(\lambda)=A_\nu(\lambda)(\lambda I_d-R)^\nu,$
%where $A_\nu(\lambda)$ is a monic $\lambda$-matrix of degree $p-\nu$ and order $d$.\\
If $A_R(R)=0_{d\times d}$ we simply say that $R$ is a \textsl{right solvent} of $A(\lambda)$.
A right solvent $R$ of $A(\lambda)$ is called \textsl{regular}
if
$\sigma(R)\cap\sigma(A^{(1)}(\cdot))=\emptyset,$
where $A^{(1)}(\lambda)$ is a  monic $\lambda$-matrix  of degree $p-1$ satisfying
$A(\lambda)=A^{(1)}(\lambda)(\lambda I_d-R)$.
\end{definition}

\begin{definition}
A set of right solvents $R_1,\ldots,R_{\mu}\in\C^{d\times d}$ of the $\lambda$-matrix $A(\lambda)$ of degree $p$ is called \textsl{ complete }
if $\sigma(A(\cdot))=\bigcup_{j=1}^{\mu}\sigma(R_j)$, where $\sigma(R_j)=\{\lambda\in\C:\det(\lambda I_d-R_j)=0\}$
is the spectrum of $R_j$. In this case,
$\nu_1+\ldots+\nu_{\mu}=p$ where $\nu_k$ is the multiplicity of the right solvent $R_k$.
\end{definition}
The Vandermonde matrix is extended in the next definition.
\begin{definition}
\label{defConfluentVandermonde matrix}
\index{Block Vandermonde matrix! Confluent Vandermonde matrix}
Suppose $R_1,\ldots,R_{\mu}$ are a complete set of right solvents of the matrix polynomial $A(\lambda)$ with multiplicities
$\nu_1,\ldots,\nu_{\mu}$, respectively. We define the  \textsl{confluent
Vandermonde matrix} $W:=W(R_1,\ldots,R_{\mu})\in \C^{pd\times pd}$ by
$W=[W_1,\ldots,W_\mu]$, where  for $k=1,\ldots,\mu$,
\begin{align*}
%\label{eqConfluentVandermonde}
&W_k=\begin{pmatrix}
I_d & 0_d &  \ldots & 0_d\\
R_k & I_d &  ~ & \vdots \\
R_k^2 & 2R_k & \ddots &  0_d \\
%R_k^3 & 3R_k^2 & 3R_k & \ddots & 0_d \\
\vdots & \vdots  & ~ & I_d \\
\vdots & \vdots  & ~ & \vdots \\
R_k^{p-1} & (p-1)R_k^{p-2} &  \ldots & \left({p-1}\atop {\nu_k-1}\right)
R_k^{p-\nu_k}
\end{pmatrix}\in\C^{pd\times\nu_kd}.
\end{align*}
\end{definition}
In the case $\mu=p$ and  $\nu_1=\ldots=\nu_p=1$ the confluent
Vandermonde matrix reduces to the classical block Vandermonde matrix
$V(R_1,\ldots,R_p)=W(R_1,\ldots,R_p)$.
%\begin{align*}
%V(R_1,\ldots,R_p):=W(R_1,\ldots,R_p)=\begin{pmatrix}
%I_d & I_d &  \ldots & I_d\\
%R_1 & R_2 & \ldots & R_p \\
%\vdots & \vdots & ~ & \vdots\\
%R_1^{p-1} & R_2^{p-1}  & \ldots & R_p^{p-1}
%\end{pmatrix}\in \C^{pd\times pd}.
%\end{align*}

%Hence, a complete set of regular right solvents $\{R_k :k=1,\ldots,p\}$ of
%$A(\lambda)$ is given by
%\begin{align*}
%\sigma(R_k)\cap\sigma(R_j)=\emptyset,~ k\not=
%j,~j,k=1,\ldots,p
%\quad\text{and}\quad
%\sigma(A(\cdot))=\bigcup_{k=1}^p\sigma(R_k).
%\end{align*}
%A criterion for a complete set of right solvents is also available:

\begin{lemma}[Maroulas \cite{Maroulas1985}, Theorem 3.4] \label{lemma:2.8}
Let $R_1,\ldots,R_{\mu}$ be right solvents of a monic $\lambda$-matrix $A(\lambda)$ of multiplicities
$\nu_1,\ldots,\nu_{\mu}$, respectively. Then $W(R_1,\ldots,R_{\mu})$ is non-singular if and only if
\begin{align*}
\sigma(A(\cdot))=\bigcup_{j=1}^{\mu}\sigma(R_j)\quad \text{ and }\quad
\sigma(R_j)\cap\sigma(R_i)=\emptyset \quad \text{ for } j,i=1,\ldots,\mu,~j\not=i.
\end{align*}
%In particular, $R_1,\ldots,R_{\mu}$ is a complete set of regular right solvents.
% In this case the
%monic $\lambda$-matrix $A(\lambda)$ has a complete set of right solvents.
\end{lemma}
Thus, we have the following relation between the solvents
of the $\lambda$-matrix $A(\lambda)$ and the coefficient matrices $A_1,\ldots,A_p$
of $A(\lambda)$.

\begin{lemma}[Maroulas \cite{Maroulas1985}] \label{lemma:2.7}
Let $R_1,\ldots,R_p$ be a complete set of regular right solvents of the monic $\lambda$-matrix $A(\lambda)=I_d\lambda^p+A_1\lambda^{p-1}+\ldots+A_{p-1}\lambda+A_p$. Then
\beao
    [A_p,\ldots,A_1]&=&-[R_1^p,\ldots,R_p^p]V^{-1}(R_1,\ldots,R_p) \quad\quad \text{ and }\\
    A(\lambda)&=&(\lambda I_d-R_p^*)\cdots (\lambda I_d-R_{2}^*)(\lambda I_d-R_1),
\eeao
where for $k=2,\ldots,p$,
\beao
    R_k^*=M_k(R_k)R_kM_k^{-1}(R_k) \quad \text{ and } \quad
    M_k(R_k)=(\lambda I_d-R_{k-1})\cdots(\lambda I_d-R_{1}).
\eeao
\end{lemma}
Interesting is that in the multivariate setting $R_k^*$
is not necessarily equal to $R_k$ for $k=1,\ldots,p$, as in
the one-dimensional case $d=1$.

\begin{definition}
\label{defRationallambdaMatrixFunc}
%$\mbox{}$
%\begin{enumerate}[(a)]
%\item
A \textsl{strictly proper rational left $\lambda$-matrix} $F(\lambda)$ with degree $p$ and  order
$(d,m)$ has
the representation
\begin{align*}
F(\lambda)&=A(\lambda)^{-1}B(\lambda),
\end{align*}
where
$
A(\lambda)$ is a monic $\lambda$-matrix of degree $p$ and order $d$,
and
$B(\lambda)$ is a $\lambda$-matrix of degree $p-1$ and order $(d,m)$.
 The rational $\lambda$-matrix
$F(\lambda)$ is called \textsl{irreducible} if $A(\lambda)$ and $B(\lambda)$ are left coprime.
%For an irreducible $F(\lambda)$ the roots of $\det(A(\lambda))$ are referred to as the \textsl{poles} of $F(\lambda)$. If $A(\lambda)$ has a
%complete set of regular right solvents $R_i$ for $i=1,\ldots,n$, then the right solvents $R_i$ are called the
%\textsl{regular left block poles} of $F(\lambda)$.
%A rational right $\lambda$-matrix $F(\lambda)=B_r(\lambda)A_r(\lambda)^{-1}$ is defined
%analogously. We denote by $\adj(A)$ the adjugate of $A$. An alternative representation for a rational $\lambda$-matrix $F(\lambda)$ is
%given by
%\begin{align*}
%F(\lambda)=\frac{1}{\det(A_l(\lambda))}\adj(A_l(\lambda))B_l(\lambda).
%\end{align*}
% and $\det(A_l(\lambda))=\det(A_r(\lambda))$ holds.
%Furthermore, we have
%$\adj(A_l(\lambda))B_l(\lambda)=B_r(\lambda)\adj(A_r(\lambda))$.
   If $F(\lambda)$ is irreducible and $R$ is a
\textsl{regular}  right solvent of $A(\lambda)$ then
 the \textsl{residue of the rational $\lambda$-matrix} $F(\lambda)$ at $R$ is defined by
\begin{equation*}
%\label{eqMatrixRes}
\res[F,R]:=\frac{1}{2\pi i}\oint_{\Gamma_R}F(\lambda)\,d\lambda,
\end{equation*}
where $\Gamma_R$ is a simple closed contour such that $\sigma(R)$ is contained in the interior of $\Gamma_R$
and $\sigma(A(\cdot))\setminus\sigma(R)$ is contained in the exterior of  $\Gamma_R$.
%\end{enumerate}
\end{definition}
%For a more detailed theory on matrix residues of rational  $\lambda$-matrices see, e.g.,
%Tsay and Shieh \cite{TsayShieh}.
The next result characterizes a rational left matrix function.
However, although Tsay and Shieh \cite{TsayShieh} assume that $d=m$, it is straightforward to extend the result
to the case $d\not=m$
(cf. Levya-Ramos \cite{Levya-Ramos1991}).

\begin{theorem}[Tsay and Shieh \cite{TsayShieh},  Theorem 4.1]
\label{thmGeneralMatrixResidue}\label{thmBlockPartialFractionExpansion}
Let $F(\lambda)=A(\lambda)^{-1}B(\lambda)$ be a irreducible strictly proper rational left $\lambda$-matrix
of degree $p$ and order $(d,m)$, and $A(\lambda)$
%has a set of regular right solvents
%$\{R_k:k=1,\ldots,\mu\}$, $\mu\leq p$, and $\bigcup_{k=1}^{\mu}\sigma(R_k)$ lies in the interior of  a simple closed contour
%$\Gamma$.
%\begin{enumerate}[(a)]
%\item  Define $F^{(k)}(\lambda):=(\lambda I_d-R_k)F(\lambda)$.
%Then
%\begin{align*}
%\label{eqGeneralMatrixResidue}
%&\frac{1}{2\pi i}\oint_\Gamma F(\lambda)\,\dif\lambda=\sum_{k=1}^{\mu} \res[F,R_k]=\sum_{k=1}^{\mu}
%F^{(k)}(R_k).
%\end{align*}
%\item
%If $F(\lambda)$ is further irreducible and $A(\lambda)$
 has a
complete set of regular right solvents $\{R_k:k=1,\ldots,p\}$. Then
\begin{align*}
%\label{eqBlockPartialFractionExpansion}
&F(\lambda)=\sum_{k=1}^p (\lambda I_d-R_k)^{-1}\res[F,R_k].
\end{align*}
%\end{enumerate}
\end{theorem}
%Finally, we recall another theorem of Tsay and Shieh \cite{TsayShieh}, which is going to be the key result for the proof of the decomposition presented in the next section. It enables us to separate a rational $\lambda$-matrix
% into a sum using the right solvents and the corresponding matrix residual.
% \begin{theorem}[Tsay and Shieh (1982), Theorem 4.1]
%\label{thmBlockPartialFractionExpansion}
%If $F(\lambda)$ is a strictly proper, irreducible rational $\lambda$-matrix, where $A(\lambda)$ has a
%complete set of regular right solvents $R_1,\ldots R_n$, then
%\begin{align}
%\label{eqBlockPartialFractionExpansion}
%&F(\lambda)=\sum_{k=1}^n (\lambda I_m-R_k)^{-1}\res[F,R_k].
%\end{align}
%\end{theorem}
%

An assumption of \Cref{thmBlockPartialFractionExpansion} is that the right solvents are regular which
excludes right solvents with multiplicities. %However, in this case a representation of $F(\lambda)$
%is derived in Shieh et al. \cite{Shieh1986}.

A formula for the calculation of a matrix residue is given in Levya-Ramos \cite[Section 6, eq. (6.13)]{Levya-Ramos1991}: Suppose
the  strictly proper left $\lambda$-matrix $F(\lambda)=A(\lambda)^{-1}B(\lambda)$ is irreducible and $A(\lambda)$
has a complete set of  regular right solvents $\{R_k:k=1,\ldots,p\}$.
Notice, the matrix $A^{\#}$ as defined in \eqref{A hash} is non-singular because
$A^{\#}$ has the only eigenvalue 1.
%Define
%the block Vandermonde matrix
Then due to \Cref{lemma:2.8} the Vandermonde matrix $V(R_1,\ldots,R_p)$ is non-singular (cf. Levya-Ramos \cite[Definition 4]{Levya-Ramos1991})
and
\begin{align} \label{Res}
\begin{pmatrix}
\res[F,R_1] \\
%\res[F,R_2]  \\
\vdots \\
\res[F,R_p]
\end{pmatrix}
=V(R_1,\ldots,R_p)^{-1}
[A^{\#}]^{-1}B^{\#}.
%\begin{pmatrix}
%I_d  & \ldots & ~ & 0_d\\
%A_1  &\ddots & ~ & \vdots \\
%\vdots  & \ddots & \ddots & \vdots \\
%%A_{p-2} & ~ & \ddots & I_d & 0_d  \\
%A_{p-1}  & \ldots & A_1 & I_d
%\end{pmatrix}^{-1}
%\begin{pmatrix}
%B_0\\
%B_2  \\
%\vdots \\
%B_{p-1}
%\end{pmatrix}.
\end{align}
%
%and for repeated solvents the formula changes to
%\begin{align}
%\label{eqMatrixResidueFormulaRepeatedSolvents}
%\begin{pmatrix}
%F_{1,1} \\
%\vdots\\
%F_{1,\nu_1}  \\
%\vdots \\
%F_{\mu,1}\\
%\vdots \\
%F_{\mu,\nu_\mu}
%\end{pmatrix}
%=W(R_1,\ldots,R_\mu)^{-1}
%\begin{pmatrix}
%I_m & 0_m & \ldots & ~ & 0_m\\
%A_1 & I_m &\ddots & ~ & \vdots \\
%\vdots & \ddots & \ddots & \ddots & \vdots \\
%A_{n-2} & ~ & \ddots & I_m & 0_m  \\
%A_{n-1} & A_{n-2} & \ldots & A_1 & I_m
%\end{pmatrix}^{-1}
%\begin{pmatrix}
%B_0\\
%B_2  \\
%\vdots \\
%B_{n-1}
%\end{pmatrix}.
%\end{align}
%

Finally, the question arises how to calculate the right solvents of  the $\lambda$-matrix $A(\lambda)$.
A possibility to characterize a right solvent is by right latent roots and latent vectors as is done
in Dennis et al. \cite{Dennis1976}.

\begin{definition}
\label{defLatentRootVector} Let $A(\lambda)$ be a $\lambda$-matrix of order $d$.
If $\lambda_i\in\C$ satisfies
$
\det(A(\lambda_i))=0$,
then $\lambda_i$ is called \textsl{latent root} of $A(\lambda)$. A  vector
$p_i\in\C^d$ satisfying
$
%\label{eqRightLatentVector}
A(\lambda_i)p_i=0_{d}
$
is called \textsl{right latent vector} of $A(\lambda)$ associated to the latent root $\lambda_i$.
%Similarly, $q_i$ is a \textsl{left latent vector} if
%$$q_i^\mathsf{T}A(\lambda_i)=0_{1\times m}.$$
\end{definition}

%\begin{theorem}[Dennis et al. \cite{Dennis1976}, Theorem 4.1]
%\label{LemSolvent}
%If $A(\lambda)$ has $d$ linearly independent right latent vectors $p_1,\ldots, p_d$
%corresponding to the latent roots $\lambda_1,\ldots,\lambda_d$.
%Define $P:=(p_1,\ldots,p_d)\in \C^{d\times d}$ and
%$\Lambda:=diag(\lambda_1,\ldots,\lambda_d)$.
%Then $R:=P\Lambda P^{-1}$ is a
%right solvent of $A(\lambda)$.
%\end{theorem}

%We want to summarize the results of \Cref{thmCompletesetofSolvents}, \Cref{thmGeneralMatrixResidue} and \Cref{LemSolvent}:

\begin{theorem} \label{Corollary:2.13}
    Suppose the monic $\lambda$-matrix $A(\lambda)$ has distinct latent roots $\lambda_1,\ldots,\lambda_{pd}$
     with corresponding %independent
     right latent vectors $p_1,\ldots,p_{pd}$, respectively.
    Define $P_k:=(p_{(k-1)d+1},\ldots,p_{kd})\in \C^{d\times d}$ and
$\Lambda_k:=diag(\lambda_{(k-1)d+1},\ldots,\lambda_{kd})$ for $k=1,\ldots,p$.
\begin{itemize}
    \item[(a)] Then $R_k:=P_k\Lambda_k P_k^{-1}$ for $k=1,\ldots,p$ is a complete set of regular
right solvents of $A(\lambda)$.
%    \item[(b)] There exist matrices  $S_1,\ldots,S_p$ in $\C^{d\times d}$, where $S_k$ is similar to $R_k$ for $k=1,\ldots,p$,
 %       such that $A(\lambda)=\prod_{k=1}^p(\lambda I_d-S_k)$. $S_k$ is not necessarily equal to $R_k$ for $k=1,\ldots,p$.
    \item[(b)] Suppose the  strictly proper left $\lambda$-matrix $F(\lambda)=A(\lambda)^{-1}B(\lambda)$ is irreducible, then  the residue of $F(\lambda)$  can be calculated as in  \eqref{Res} and
        \begin{align*}
            %\label{eqBlockPartialFractionExpansion}
        &F(\lambda)=\sum_{k=1}^p (\lambda I_d-R_k)^{-1}\res[F,R_k].
\end{align*}
\end{itemize}
\end{theorem}
\begin{proof}
    (a) \, is proven in
            Dennis et al. \cite{Dennis1976}, Theorem 4.5. %  and (b) in \Cref{lemma:2.7} (see \cite{Dennis1976}, Corollary 4.1).\
    (b) \, follows from (a) and \Cref{thmGeneralMatrixResidue}.
\end{proof}

\section{Results}
\label{ch:chapter3:sec:DecompMCARMA}
%% ==============================

%\begin{proposition}[Schlemm and Stelzer (2012), Corollary 3.4]
%The classes of MCARMA and causal continuous-time state space models are equivalent.
%\end{proposition}

%
%\marginpar{andere Dimensionen in \Cref{definition CARMA}}
%
%\marginpar{Wenn $N$ sonst nicht mehr verwendet wird, würde ich direkt $pm$ schreiben}
%
%
%\begin{assumptionletter}
%\label{AssDisEV}
%The eigenvalues $\lambda_1,\ldots,\lambda_N$ of $\mathcal{A}$ as in \Cref{defLSSM} and consequently of $A\in M_N(\C)$ in \Cref{definition CARMA}, are distinct, where the dimensions satisfies $N=pm$.
%\end{assumptionletter}

In this section we present criteria for a MCARMA process to be a sum
of multivariate Ornstein-Uhlenbeck processes.
For the rest of the paper we will assume the following:

\begin{assumptionletter} \label{assumption:B}
\label{AssNegEV}
Let $A(\lambda)$, $B(\lambda)$ be defined as in \eqref{Pol}
 and   $F(\lambda)=A(\lambda)^{-1}B(\lambda)$ be irreducible.
%Suppose  $\sigma(A(\cdot))\subset\{(-\infty,0)+i\R\}$ and   $\E[\log(\max(1,\|L(1)\|))]<\infty$.
%The MCARMA process $Y$ has the representation as given in \eqref{MCARMA} and hence, is stationary.
Assume further that $A(\lambda)$ has a
    complete set of regular right solvents $\{R_k:k=1,\ldots,p\}$.
\end{assumptionletter}
Instead of assuming that the   right solvents $\{R_k:k=1,\ldots,p\}$ are  complete and regular,
it is equivalent to assume that $V(R_1,\ldots,R_p)$ is non-singular (see \Cref{lemma:2.8}).
A sufficient condition for $A(\lambda)$ to have a complete set of  regular right solvents
is that $A^*$ as defined in \eqref{def:A:star} has distinct eigenvalues because
 $\sigma(A^*):=\{\lambda\in\C:\det(A^*-\lambda I_{pd})=0\}=\sigma(A(\cdot))$, due
to Marquardt and Stelzer~\cite[Lemma 3.8]{MarquardtStelzer2007},
such that by \Cref{Corollary:2.13} %to \autoref{AssNegEV} we can apply Dennis et al. \cite[Theorem 4.1]{Dennis1976}  and hence,
  the statement follows. However, this is only a sufficient but not a necessary assumption.

\begin{theorem}
\label{thmBlockDecompStateSpaceMod}
%Suppose that $A^*$ as defined in \eqref{def:A:star} has distinct eigenvalues $\lambda_1,\ldots,\lambda_{pd}$.
%\begin{enumerate}[(a)]
%    \item The autoregressive matrix polynomial $A(\lambda)$ as given in \eqref{PolP} has a complete regular set of right solvents $(R_k)_{k=1,\ldots,p}$.
%    \item Suppose $B(\lambda)$ is the moving-average matrix polynomial as given in \eqref{PolQ} and $F(\lambda):=A(\lambda)^{-1}B(\lambda)$ is irreducible.
     % The   process $Y$  has the representation as
%      The process , $t\geq 0$, as sum of complex-valued possible dependent multivariate Ornstein-Uhlenbeck processes $Y_1,\ldots,Y_p$  where
Define for $k=1,\ldots,p$ the multivariate complex-valued Ornstein-Uhlenbeck processes
\begin{align}
%\intertext{where for $s,t\in \R$ with $s<t$}
\label{eqDecompStateSpaceModBlockcomponent}
Y_k(t)&=\e^{R_k t}Y_k(0)+\int_0^t\e^{R_k(t-u)}\res[F,R_k]\,\dif L(u), \quad t\geq 0,
\end{align}
with some initial condition $Y_k(0)$ in $\C^d$ such that $V(R_1,\ldots,R_p)[Y_1(0)^\top,\ldots,Y_p(0)^\top]^\top\in\R^{pd}$. Then $Y(t)=\sum^p_{k=1}Y_k(t)$
is an $\R^d$-valued solution of the state space model  \eqref{state:space} and hence, a MCARMA$(p,q)$-process.
\end{theorem}

\begin{proof}
Of course,
\beao
    Y(t)=C^*\e^{A^*t}X(0)+\int_0^tC^*\e^{A^* (t-u)}B^*\,\dif L(u)
\eeao
is an $\R^d$-valued solution of the state space model~\eqref{state:space} with some initial condition $X(0)\in \R^{pd}$.
Define
\beao
    E^*:=[I_d,\ldots,I_d]\in\R^{d\times pd}, \quad
    F^*=\left(\begin{array}{c}\res[F,R_1] \\ \vdots \\\res[F,R_p]
    \end{array}\right)\in\C^{pd\times d} \,
    \text{ and }\,R^*:=\diag(R_1,\ldots,R_p)\in\C^{pd\times pd}
\eeao
 as a block diagonal matrix.
%Due to \Cref{aux_Lemma}(b) we receive
%\beam \label{eq2}
%    Y(t)=C^*\e^{A^*t}X(0)+\sum_{k=1}^p\int_0^t\e^{R_k(t-u)}\res[F,R_k]\,\dif L(u).
%\eeam
%Next, we distinguish the different assumptions in \autoref{AssumptionB1} and show that
%$Y(t)=\sum_{k=1}^pY_k(t)$.
%\begin{enumerate}
%    \item[(B1)] Suppose $Y_1(0)=\ldots=Y_p(0)=0_d$. If we take $X(0):=0_{pd}$ then
%        \beao
%         Y(t)=\sum_{k=1}^p\int_0^t\e^{R_k(t-u)}\res[F,R_k]\,\dif L(u)=\sum_{k=1}^pY_k(t), \quad t\geq 0.
%        \eeao
 %   \item[(B2)] % From \Cref{aux_Lemma} we already know that $ E^*(\lambda I_{pd}-R^*)^{-1}F^*
        %=A(\lambda)^{-1}B(\lambda)
     %   =C^*(\lambda I_{pd} -A^*)^{-1}B^*$.   A consequence of
     %   Bernstein~\cite[Proposition 12.9.8]{Bernstein2009}
     %   (due to \cite[Corollary 12.9.15 and Theorem 12.9.16]{Bernstein2009} and the irreducibility
      %  of $F(\lambda)$ the assumptions are satisfied) % Hannan~\cite[Theorem 2.3.4]{hannandeistler2012}
%   is that there exists a non-singular matrix $T\in \R^{pd\times pd}$
%  \textcolor{red}{komplexwertig??? sonst muss man $R_k$ reellwertig annehmen}
Due to  \eqref{Res} and \eqref{2} the relation
$$F^*=V(R_1,\ldots,R_p)^{-1}[A^{\#}]^{-1}B^{\#}=V(R_1,\ldots,R_p)^{-1}B^*$$ holds. A further inspection of the matrices give
\beao \label{T}
\begin{array}{rcl}
    A^*V(R_1,\ldots,R_p)=V(R_1,\ldots,R_p)R^* \quad \text{ and } \quad C^*V(R_1,\ldots,R_p)=E^*, %\quad \text{ and } \quad\\
    %V(R_1,\ldots,R_p)F^*&=&.
\end{array}
\eeao
where we used that $R_k$ is a right solvent of $A(\lambda)$.
Therefore, define $T:=V(R_1,\ldots,R_p)$, $Y^*(0):=[Y_1(0)^\top,\ldots,Y_p(0)^\top]^\top$ and
$X^*(0):=TY^*(0)\in\R^{pd}$ such that
  \beao
    A^*=TR^*T^{-1}, \quad B^*=TF^* \quad \text{ and } \quad C^*=E^*T^{-1}.
  \eeao
  In particular, $\e^{A^*t}=T\e^{R^*t}T^{-1}$, $t\in\R$. Then
  for $t\geq 0$,
  \beao
    Y(t)&=&C^*\e^{A^* t}X^*(0)
        +\int_0^tC^*\e^{A^* (t-u)}B^*\,\dif L(u)\\
        &=&
        E^*\e^{R^* t}T^{-1}TY^*(0)
        +\int_0^tE^*\e^{R^* (t-u)}F^*\,\dif L(u)\\
        &=&\sum_{k=1}^p\left[\e^{R_k t}Y_k(0)+\int_0^t\e^{R_k (t-u)}\res[F,R_k]\,\dif L(u)\right]\\
        &=&\sum_{k=1}^pY_k(t)
  \eeao
  is $\R^d$-valued.
%      \item[(B3)] Suppose $Y_k(0):=\int_{-\infty}^0\e^{-R_ku}\,\dif L(u)$ for $k=1,\ldots,p$. Define $X(0):=\int_{-\infty}^0\e^{-A^*u}B^*\dif L(u)$, which is well-defined due to Sato and Yamazato~\cite[Theorem 4.1]{Sato:Yamazato:1984}. Then
%    an application of \Cref{aux_Lemma}(b) gives
%        \beao
%            C^*\e^{A^*t}X(0)=\int_{-\infty}^0C^*\e^{A^*(t-u)}B^*\dif L(u)
%            =\sum_{k=1}^p\int_{-\infty}^0\e^{R_k (t-u)}\res[F,R_k]\,\dif L(u)=\sum_{k=1}^p\e^{R_k t}Y_k(0).
%        \eeao
%    Plugging this in \eqref{eq2} gives $Y(t)=\sum_{k=1}^pY_k(t)$, $t\geq 0$.
%\end{enumerate}
%It is straightforward to see that \eqref{eqDecompStateSpaceModBlockcomponent} is satisfied in all
%three cases.
\end{proof}

\begin{remark} $\mbox{}$
\begin{itemize}
\item[(a)] If $\sigma(A^*)$ has only distinct eigenvalues then \Cref{Corollary:2.13}
    gives the possibility to calculate a complete set of regular right solvents. Due to \Cref{Res}
    we are able to calculate the residues as well. Thus, we obtain via \eqref{eqDecompStateSpaceModBlockcomponent}
    a representation of the MCARMA process as sum of Ornstein-Uhlenbeck processes.
%\item[(a)] In (B3), the assumption that $\sigma(A(\cdot))\subset\{(-\infty,0)+i\R\}$  is equivalent to
%the eigenvalues of $A^*$ having strictly negative real parts. In this case, $Y$ as defined in
%\Cref{thmBlockDecompStateSpaceMod} is a stationary MCARMA process (cf. Marquardt and Stelzer
%\cite{MarquardtStelzer2007}).
\item[(b)] Since the solvents $R_1,\ldots,R_p$ are not unique, the representation of $Y$ as sum of
    Ornstein-Uhlenbeck processes is not unique as well (cf. \Cref{Example}), only in the case $d=1$ we have uniqueness.
\item[(c)] Any linear combination $\sum_{k=1}^pL_kY_k(t)$, $t\geq 0$, where $L_1,\ldots,L_p\in\R^{d\times d}$, of $\R^d$-valued multivariate Ornstein-Uhlenbeck processes
     $Y_1,\ldots,Y_p$ is  a MCARMA$(p,p-1)$-process.
    But the exponent $R_k$ in the definition of $Y_k$ is not necessarily a right solvent of the autoregressive polynomial
    of the MCARMA process.
    This is essential to derive a VARMA representation of the discrete-time sampled MCARMA process
    later on.
%\item[(e)] Unfortunately, the result of   Horn and Johnson~\cite[Proposition 12.9.8]{horn1994topics}, which can
%    be found in many standard textbooks on linear systems and control theory, is only valid for
%    real-valued matrices and not for complex-valued matrices.
\end{itemize}
\end{remark}

\begin{corollary} \label{corollary}
Suppose  $\sigma(A(\cdot))\subset\{(-\infty,0)+i\R\}$ and   $\E[\log(\max(1,\|L(1)\|))]<\infty$.
Define for $k=1,\ldots,p$ the multivariate complex-valued Ornstein-Uhlenbeck processes
\begin{align*}
%\intertext{where for $s,t\in \R$ with $s<t$}
Y_k(t)&=\int_{-\infty}^t\e^{R_k(t-u)}\res[F,R_k]\,\dif L(u), \quad t\in\R,
\end{align*}
 Then $Y(t)=\sum^p_{k=1}Y_k(t)=\int_{-\infty}^t\sum_{k=1}^p\e^{R_k(t-u)}\res[F,R_k]\,\dif L(u)$, $t\in\R$,
is a stationary $\R^d$-valued solution of the state space model  \eqref{state:space} and hence, a MCARMA$(p,q)$-process.
\end{corollary}
Due to Sato and Yamazato~\cite[Theorem 4.1]{Sato:Yamazato:1984} the stationary Ornstein-Uhlenbeck processes $Y_k$ are well-defined.

\begin{remark} $\mbox{}$
\begin{enumerate}
\item[(a)] Let $\Gamma_k$ be a simple closed contour such that $\sigma(R_k)$ lies in the interior of $\Gamma_k$ and the residuary spectrum $\sigma(A(\cdot))\setminus\sigma(R_k)$ lies in the exterior of $\Gamma_k$ and $\Gamma:=\bigcup_{k=1}^p\Gamma_k$.
Due to Cauchy's integral formula (see Lax \cite[Theorem 17.5]{lax2002}), and \Cref{thmBlockPartialFractionExpansion} we obtain for $t\geq 0$,
\begin{eqnarray*}
\sum_{k=1}^p \e^{tR_k}\res[F,R_k]=\frac{1}{2\pi i}\sum_{k=1}^p\oint_{\Gamma} \e^{t\lambda}(\lambda I_d-R_k)^{-1}\res[F,R_k]\dif \lambda
  %  &=&\frac{1}{2\pi i}\oint_{\Gamma} \e^{t\lambda}E^*(\lambda I_{pd}-R^*)^{-1}F^*\,\dif\lambda\\
    =\frac{1}{2\pi i}\oint_{\Gamma} \e^{t\lambda}F(\lambda)\,\dif\lambda.
\end{eqnarray*}
In particular, if  $\sigma(A(\cdot))\subset\{(-\infty,0)+i\R\}$ then the kernel function satisfies for $t\geq0$,
\beao
    \sum_{k=1}^p \e^{tR_k}\res[F,R_k]&=&\frac{1}{2\pi }\int_{-\infty}^\infty \e^{ti\omega}F(i\omega)\,\dif\omega.
\eeao
\item[(b)] In the case of repeated right solvents, $Y$ has  not the representation as a sum of multivariate
Ornstein Uhlenbeck processes due to the representation of $F(\lambda)=\sum_{k=1}^\mu\sum_{j=1}^{\nu_k}(\lambda
I_d-R_k)^{-j}F_{k,j}$
 in Shieh et. al. \cite{Shieh1986} % \Cref{thmBlockPartialExpRepeatedBlockPoles}
 and hence,
\beao
    \sum_{k=1}^p \e^{tR_k}\res[F,R_k]&\not=&\frac{1}{2\pi i }\oint_{\Gamma} \e^{t\lambda}F(\lambda)\,\dif\lambda
\eeao
(cf.  Brockwell and Lindner~\cite[Lemma 2.4]{Brockwell2009} in the case of one-dimensional
CARMA processes).
\end{enumerate}
\end{remark}

\begin{comment}
\begin{corollary}
Let $R_1,\ldots,R_p\in\R^{d\times s}$ be a complete set of regular right solvents of the monic $\lambda$-matrix $A(\lambda)$
and $F_1^*,\ldots,F_p^*\in\R^{d\times m}$. Suppose for $k=1,\ldots,p$ that $Y_k$ is a multivariate
Ornstein-Uhlenbeck process of the form
\begin{align*}
    Y_k(t)&=\e^{R_k t}Y_k(0)+\int_0^t\e^{R_k(t-u)}F_k^*\,\dif L(u), \quad t\geq 0
\end{align*}
with $Y_1(0),\ldots,Y_p(0)\in\R^{d}$.
Then $Y(t)=\sum^p_{k=1}Y_k(t)$ is a MCARMA$(p,p-1)$ process with autoregressive polynomial
$A(\lambda)$ and moving average polynomial $B(\lambda)=B_0\lambda^{p-1}+\ldots+ B_{p-2}\lambda+B_{p-1}$ where
\begin{align*}
\left(\begin{array}{c}
    B_0\\
    \vdots\\
    B_{p-1}
    \end{array}\right)
= A^{\#}V(R_1,\ldots,R_p)F^*
\quad \text{ and } \quad
F^*=
\begin{pmatrix}
F_1^* \\
%\res[F,R_2]  \\
\vdots \\
F_p^*
\end{pmatrix}.
\end{align*}
\end{corollary}
\begin{proof}
  Notice that $F_k^*=\res[F,R_k]$ for $k=1,\ldots,p$ due to \eqref{Res} and that
  all matrices have real-valued entries. The rest follows from \Cref{thmBlockDecompStateSpaceMod}.
\end{proof}

\begin{remark} $\mbox{}$
\begin{itemize}
    \item[(a)] The assumption that $R_1,\ldots,R_p$ are real-valued matrices is sufficient that
    $B_0,\ldots,B_{p-1}$ are real-valued matrices. Otherwise they might be complex-valued and $Y(t)$
    might be $\C^d$-valued as well.
    \item[(b)]
\end{itemize}
\end{remark}
\end{comment}

\begin{example} \label{Example}
Let
\beao
    A(\lambda)=I_d\lambda^2+\left(\begin{matrix} -11 & 22\\ -12 &  21 \end{matrix}\right)\lambda
            +\left(\begin{matrix}-42 & 52 \\ -36 & 44\end{matrix}\right)
            \quad \text{ and } \quad
    B(\lambda)=I_d, \quad \lambda\in\C,
\eeao
be given. Then
\beao
    R_1= \left(\begin{matrix} 0 & -1\\ 2 & -3\end{matrix}\right), \quad\quad
    R_2= \left(\begin{matrix} -3 & -2 \\ 0 & -4\end{matrix}\right),\quad\quad
    R_3= \left(\begin{matrix} -7 & 6\\ -3 & 2 \end{matrix}\right), \quad\quad
    R_4= \left(\begin{matrix}  -3 & 0.5\\ 0 & -2\end{matrix}\right)
\eeao
are right solvents of $A(\lambda)$. The pair $R_1, R_2$ and the pair $R_3,R_4$, respectively
build a complete set of regular right solvents of $A(\lambda)$. Then \Cref{thmBlockDecompStateSpaceMod} and
the formula for the residues \eqref{Res}
give that both $Y_1(t)+Y_2(t)$ with
\beao
    Y_1(t)&=&\e^{R_1 t}Y_1(0)+\int_0^t\e^{R_1(t-u)}\left(\begin{matrix} 1 & -1\\ -2 & 3\end{matrix}\right)\,\dif L(u),\\
    Y_2(t)&=&\e^{R_2 t}Y_2(0)+\int_0^t\e^{R_2(t-u)}\left(\begin{matrix} -1 & 1\\ 2 & -3\end{matrix}\right)\,\dif L(u),
\eeao
and $Y_3(t)+Y_4(t)$ with
\beao
   Y_3(t)&=&\e^{R_3 t}Y_3(0)+\int_0^t\e^{R_3(t-u)}\left(\begin{matrix} 8 & -11\\ 6 & -8\end{matrix}\right)\,\dif L(u),\\
   Y_4(t)&=&\e^{R_4 t}Y_4(0)+\int_0^t\e^{R_4(t-u)}\left(\begin{matrix}-8 & 11\\ -6 & 8\end{matrix}\right)\,\dif L(u),
\eeao
are MCARMA$(2,0)$ processes with AR polynomial $A(\lambda)$ and MA polynomial $B(\lambda)$.
\end{example}

For the rest of the paper we assume:

\begin{assumptionletter}
\label{Ass:C}
 $Y$ has the representation as given in \Cref{thmBlockDecompStateSpaceMod}, $\E\|L(1)\|^2<\infty$ and $\E L(1)=0_m$.
\end{assumptionletter}
Now, we are able  to present an
alternative representation of the covariance function of a stationary MCARMA process. % which depends on the residuals of the rational matrix function $F(\lambda)$.
\begin{proposition} \label{Propo_acf}
%Suppose %that $A^*$ has distinct eigenvalues $\lambda_1,\ldots,\lambda_{pd}$.
%Furthermore, let
%$\E\Vert L(1)\Vert^2<\infty$ and $\E(L(1))=0_m$.
Suppose the setting of \Cref{corollary}. The covariance function $(\gamma_Y(l))_{l\in\Z}=(\Cov(Y(t+l),Y(t)))_{l\in\N_0}$ of $Y$ has the representation
\label{propCovSumOU}
\begin{eqnarray*}
\gamma_Y(l)=\sum_{i=1}^p \e^{lR_i}\Sigma_i, \quad l\geq 0, \quad \text{ where } \quad \Sigma_i:=\sum_{j=1}^p\int_0^\infty\e^{uR_i}\res[F,R_i]\Sigma_L\res[F,R_j]^\mathsf{H} \e^{uR_j^\mathsf{H}}\,\dif u
%\label{eqCovSumOU}
\end{eqnarray*}
and for a matrix $Z\in\C^{d\times d}$ we denote by  $Z^{\mathsf{H}}$ the transposed complex conjugated of $Z$.
\end{proposition}
\begin{proof}
An application of \Cref{corollary} gives
\begin{align*}
\gamma_Y(l)%&%\Cov(Y(t+l),Y(t))
%\\
%\Cov\left(\sum_{i=1}^p Y_i(t+l),\sum_{j=1}^pY_j(t)\right)
%\\
=&\sum_{i,j=1}^p\Cov\left(\int_{-\infty}^{t+l}\e^{R_i(t+l-u)}\res[F,R_i]\,\dif L(u), \int_{-\infty}^t\e^{R_j(t-u)}\res[F,R_j]\,\dif L(u)\right)
\\
=&\sum_{i,j=1}^p \e^{lR_i} \int_0^\infty\e^{uR_i}\res[F,R_i]\Sigma_L\res[F,R_j]^\mathsf{H} \e^{uR_j^\mathsf{H}}\,\dif u,
\end{align*}
which completes the proof.
\end{proof}

A final aim  is to derive a VARMA representation for a MCARMA process observed at discrete time-points.
% where we use the representation of $Y$ as a sum of multivariate Ornstein-Uhlenbeck processes as given in \Cref{thmBlockDecompStateSpaceMod}.
%We use the representation
%  $Y(t)=\sum^p_{k=1}Y_k(t)$ where
%\begin{align}
%\intertext{where for $s,t\in \R$ with $s<t$}
%\label{eqDecompStateSpaceModBlockcomponent}
%Y_k(t)&=\int_{-\infty}^t\e^{R_k(t-u)}\res[F,R_k]\,\dif L(u), \quad t\in\R, \quad k=1,\ldots,p,
%\end{align}
% derived in \Cref{thmBlockDecompStateSpaceMod}.
 To distinguish the notation between the continuous-time process and the sampled discrete-time process, we write $Y_n^{(h)}$ for $Y(nh)$ in the following and accordingly $Y_{k,n}^{(h)}$ for $Y_k(nh)$ for some fixed $h>0$. Let us first state an auxiliary lemma.

\begin{lemma}
\label{LemInductionRepresentation}
For any $k=1,\ldots,p$, $n\geq p$, $l=0,\ldots,n$ and any  matrices $C_1,\ldots,C_l\in\C^{d\times d}$ it holds that
\begin{eqnarray*}
Y_{k,n}^{(h)}=\sum_{r=1}^lC_rY_{k,n-r}^{(h)}+\left(\e^{hlR_k}-\sum_{r=1}^lC_r\e^{h(l-r)R_k}\right)Y_{k,n-l}^{(h)}
+\sum_{r=0}^{l-1}\left(\e^{hrR_k}-\sum_{j=1}^rC_j\e^{h(r-j)R_k}\right)N_{k,n-r}^{(h)},
\end{eqnarray*}
where
$N_{k,n}^{(h)}=\int_{(n-1)h}^{nh}\e^{R_k(nh-u)}\res[F,R_k]\,\dif L(u)$.
\end{lemma}
\begin{proof}
The proof goes in the same vein as the proof of equation (2.8) in  Brockwell and Lindner~\cite{Brockwell2009} for scalars
$c_1,\ldots,c_l$ instead of matrices
$C_1,\ldots,C_l$, since $Y_k$ is a multivariate Ornstein-Uhlenbeck process.
 \begin{comment}
 We show the claim by induction. The assertion is clear for $l=0$ since we always set the empty sum to zero. Assume the equation holds for some $l\in\N$, i.e.,
 \begin{align*}
Y_{k,n}^{(h)}=&\sum_{r=1}^lC_rY_{k,n-r}^{(h)}+\left(\e^{hlR_k}-\sum_{r=1}^lC_r\e^{h(l-r)R_k}\right)Y_{k,n-l}^{(h)}
+\sum_{r=0}^{l-1}\left(\e^{hrR_k}-\sum_{j=1}^rC_j\e^{h(r-j)R_k}\right)N_{k,n-r}^{(h)}.
\end{align*}
Due to \eqref{eqDecompStateSpaceModBlockcomponent} the equality
$Y_{k,n-l}^{(h)}=\mathrm{e}^{hR_k}Y_{k,n-(l+1)}^{(h)}+N_{k,n-l}^{(h)}$ holds,
and hence,
\begin{align*}
Y_{k,n}^{(h)}=&\sum_{r=1}^lC_rY_{k,n-r}^{(h)}+\bigg[\e^{hlR_k}-\sum_{r=1}^lC_r\e^{h(l-r)R_k}\bigg]
\Big(\mathrm{e}^{hR_k}Y_{k,n-(l+1)}^{(h)}+N_{k,n-l}^{(h)}\Big)
\\
&
+\sum_{r=0}^{l-1}\bigg[\e^{hrR_k}-\sum_{j=1}^rC_j\e^{h(r-j)R_k}\bigg]N_{k,n-r}^{(h)}
\\
=&\sum_{r=1}^{l+1}C_rY_{k,n-r}^{(h)}+\bigg[\e^{h(l+1)R_k}-\sum_{r=1}^lC_r\e^{h(l+1-r)R_k}-C_{l+1}\bigg]Y_{k,n-(l+1)}^{(h)}
\\
&
+\sum_{r=0}^l\bigg[\e^{hrR_k}-\sum_{j=1}^rC_j\e^{h(r-j)R_k}\bigg]N_{k,n-r}^{(h)},
\end{align*}
which completes the induction step.
\end{comment}
\end{proof}

Eventually, we obtain a VARMA$(p,p-1)$ representation for the sampled version of a \linebreak MCARMA$(p,q)$ process.

\begin{theorem}
\label{thmBlockweakVARMArep}
%Let $S_k$ be similar matrices to the right solvents $R_k$, $k=1,\ldots p$, such that the
Define
\beao
    \Psi_0:=I_d, \quad \quad [\Psi_p,\ldots,\Psi_1]:=-[\e^{-phR_1},\ldots,\e^{-phR_p}]V^{-1}(\e^{-hR_1},\ldots,\e^{-hR_p})\in\C^{p\times pd}
\eeao
and the $\lambda$-matrix  $\Phi(\lambda):= I_d-\Phi_1\lambda-\ldots-\Phi_p\lambda^p$
of degree $p$ and order $d$
with
$
    \Phi_j:=-\Psi_p^{-1}\Psi_{p-j}  \text{ for } j=1,\ldots,p.
$
Then there exists a   $\lambda$-matrix $\Theta(\lambda)=I_d+\Theta_1\lambda+\ldots+\Theta_{p-1}\lambda^{p-1}$ of degree  $p-1$ and order $d$ such that
\begin{align}
\label{eqDiscrVARMA}
\Phi(B)Y_n^{(h)}=\Theta(B)\varepsilon_n^{(h)},\quad n\geq p,
\end{align}
where $B$ denotes the backshift operator (i.e. $B^jY_n^{(h)}=Y_{n-j}^{(h)}$ for $j\in\N$) and $(\varepsilon_n^{(h)})_{n\geq p}$ is
a $d$-dimensional weak white noise. Thus, $(Y_n^{(h)})_{n\geq p}$ admits a weak $VARMA(p,p-1)$ representation.
\end{theorem}
%Due to Maroulas \cite[Theorem 5.2.]{Maroulas1985}, we use only similar matrices $S_i$ in the factorization instead of the right solvents $R_i$.
\begin{proof}
First, we will show that $\Phi(\lambda)$ is well-defined and has the complete set of regular right solvents
$\e^{-hR_1}$, \ldots, $\e^{-hR_p}$. Due to \autoref{assumption:B} and \Cref{lemma:2.8}, the
Vandermonde matrix $V(\e^{-hR_1},\ldots,\e^{-hR_p})$ is non-singular and finally, $\Psi_1,\ldots,\Psi_p$
is well-defined. A conclusion of \autoref{assumption:B} and  \Cref{lemma:2.7} is then that $\e^{-hR_1}$, \ldots, $\e^{-hR_p}$
is a complete set of regular right solvents of  $\Psi(\lambda)=I_d\lambda^p
+\Psi_1\lambda^{p-1}+\ldots+\Psi_p$. Note that $\Psi_p=(-1)^p\e^{-h R_p^*}\cdot\ldots\cdot \e^{-h R_2^*}\cdot\e^{-h R_1}$
where $R_2^*,\ldots,R_p^*$ are defined as in \Cref{lemma:2.7}. Since the eigenvalues of $\e^{-h R_k^*}$, $k=2,\ldots,p$
and $\e^{-h R_1}$  are non-zero, the matrix $\Psi_p$ is non-singular.
Finally, $\Phi(\lambda)=\Psi_p^{-1}\Psi(\lambda)$ is well-defined and
has the complete set of regular right solvents
$\e^{-hR_1}$, \ldots, $\e^{-hR_p}$.

Due to \eqref{eqDecompStateSpaceModBlockcomponent} we obtain that $Y_n^{(h)}=\sum_{k=1}^p Y_{k,n}^{(h)}$, $n\geq p$, where %$Y_{k,n}^{(h)}:=Y_k(nh)$ and
\begin{eqnarray*}
Y_{k,n}^{(h)}=\e^{hR_k}Y_{k,n-1}^{(h)}+N_{k,n}^{(h)} \quad \text{and} \quad
N_{k,n}^{(h)}=\int_{(n-1)h}^{nh}\e^{R_k(nh-u)}\res[F,R_k]\,\dif L(u)
\end{eqnarray*}
%and hence,
%$
%(I_d-\e^{hR_k}B)Y_{k,n}^{(h)}=N_{k,n}^{(h)}.$
%
(cf. Schlemm and Stelzer~\cite[Lemma 5.2]{SchlemmStelzer2012}).
An application of \Cref{LemInductionRepresentation} with $l=p$ and $C_r=\Phi_r$ for $r=1,\ldots,p$  gives
\begin{eqnarray*}
Y_{k,n}^{(h)}=\sum_{r=1}^p\Phi_rY_{k,n-r}^{(h)}+\left[\e^{hpR_k}-\sum_{r=1}^p\Phi_r\e^{h(p-r)R_k}\right]Y_{k,n-p}^{(h)}
+\sum_{r=0}^{p-1}\left[\e^{hrR_k}-\sum_{j=1}^r\Phi_j\e^{h(r-j)R_k}\right]N_{k,n-r}^{(h)}.
\label{eq1proofBlockweakVARMArep}
\end{eqnarray*}
The fact that $\e^{-hR_k}$ is a right solvent of $\Phi(\lambda)$  implies that
 $$\e^{hpR_k}-\sum_{r=1}^p\Phi_r\e^{h(p-r)R_k}=\Phi(\e^{-hR_k})\e^{phR_k}=0_{d\times d}, \quad k=1,\ldots,p.$$
 Hence, we obtain
\begin{align} \label{eqARequBlocksampled}
\Phi(B)Y_{k,n}^{(h)}=Y_{k,n}^{(h)}-\sum_{r=1}^p\Phi_rY_{k,n-r}^{(h)} =\sum_{r=0}^{p-1}\left[\e^{hrR_k}-\sum_{j=1}^r\Phi_j\e^{h(r-j)R_k}\right]N_{k,n-r}^{(h)}=:U_n^{(h)}.%=:W_{k,n}^{(h)}.
\end{align}
Define for $r=0,\ldots,p-1$ the  iid sequence $(W_{r,n}^{(h)})_{n\in\Z}$ in $\C^{d}$ as
\begin{equation}\label{eqDiscIIDSequence}
W_{r,n}^{(h)}:=\int_{(n-1)h}^{nh}\sum_{k=1}^{p}\left[\e^{hrR_k}-\sum_{j=1}^{r}\Phi_j\e^{h(r-j)R_k}\right]\e^{R_k(nh-u)}\res[F,R_k]\dif L(u).
\end{equation}
%and
%\begin{equation} \label{eqARequBlocksampled}
%    U_n^{(h)}:=\sum_{r=0}^{p-1}W_{r,n-r}^{(h)}, \quad n\in\Z.
%\end{equation}
Summation over $k$ and rearranging leads to
\begin{equation*}
\Phi(B)Y_n^{(h)}=U_n^{(h)}=\sum_{r=0}^{p-1}W_{r,n-r}^{(h)}, \quad n\geq p.
\end{equation*}
Since $([W_{0,n}^{(h)\,\top},\ldots,W_{p-1,n}^{(h)\,\top}]^\top)_{n\in\Z}$ is a sequence of iid random vectors, the $d$-dimensional sequence $(U_n^{(h)})_{n\in\Z}$ is $(p-1)$-dependent.
Define
 \beao
    \varepsilon_n^{(h)}:=U_{n}^{(h)}-\mathcal{P}_{\mathcal{M}_{n-1}}U_{n}^{(h)}, \quad n\in\Z,
 \eeao
 where $\mathcal{P}_{\mathcal{M}_{n-1}}$ denotes the orthogonal projection on
 $\mathcal{M}_{n-1}:=\overline{\text{sp}}\{U_{j}^{(h)}:-\infty< j\leq n-1\}$ and the
 closure is taken in the Hilbert space of square integrable complex random vectors
 with inner product $(U_1,U_2)\mapsto \E({U}_1^{\mathsf{H}} U_2)$  for random vectors $U_1,U_2$ in $\C^d$.
 Then $\Theta_1,\ldots\Theta_{p-1}$ is given as the solution of the equation
 \beao
    \mathcal{P}_{\overline{\text{sp}}\{\varepsilon_{n-p+1}^{(h)},\ldots,\varepsilon_{n-1}^{(h)}\}}U_{n}^{(h)}
        =\Theta_1\varepsilon_{n-1}^{(h)}+\ldots+\Theta_{p-1}\varepsilon_{n-p+1}^{(h)}.
 \eeao
  As in the proof of Brockwell and Davis \cite[Proposition 3.2.1]{BrockwellDavis1998}
  for one-dimensional $(p-1)$-dependent processes we can follow then the statement.
\end{proof}

%\subsection{Autocovariance structure of the noise sequence}

\begin{remark} $\mbox{}$
\begin{itemize}
    \item[(a)] Characteristic is that the  $\lambda$- matrix $\Psi(\lambda)$ has the complete
            set of right solvents $\e^{-h R_1},\ldots,\e^{-h R_p}$ but due to \Cref{lemma:2.7}
            it has not necessarily the representation as $\prod_{k=1}^p(\lambda I_d- \e^{-h R_k})$.
            Thus, the $\lambda$-matrix $\Phi(\lambda)$ is not necessarily
            $\psi_p^{-1}\prod_{k=1}^p(\lambda \e^{h R_k}- I_d)$. This differs to
            the one-dimensional case where multiplication is commutative.
            However, $\Psi(\lambda)$ is the unique $\lambda$-matrix with right solvents $\e^{-h R_1},\ldots,\e^{-h R_p}$ and $\Psi(0)=I_d$.
    \item[(b)] If $\sigma(A(\cdot))\subset\{(-\infty,0)+i\R\}$ holds then
            \beao
                \sigma(\Psi(\cdot))=\bigcup_{k=1}^p\sigma(\e^{-hR_k})
                =\{\e^{-h\lambda}:\lambda\in\bigcup_{k=1}^p\sigma(R_k)\}=\{\e^{-h\lambda}:\lambda\in\sigma(A)\},
            \eeao
            is outside the closed unit disc and hence, $\Psi(\lambda)$ is Schur-stable.
\end{itemize}
\end{remark}

Finally, we state the covariance function of the series $U^{(h)}:=(U_n^{(h)})_{n\geq\Z}$ given in \eqref{eqARequBlocksampled}.
The second-order properties of the series $U^{(h)}$ are of interest for indirect
estimation as is done, e.g., in Fasen-Hartmann and Kimmig~\cite{Fasen:Kimmig:2019} for CARMA processes.
The basic idea is that the VARMA parameters of $(Y(nh))_{n\in\N}$ are estimated by standard techniques.
Taking identifiability issues into account the autoregressive parameters of the continuous-time  process
are then estimated from the autoregressive parameters of the discrete-time VARMA process. Finally,
 a comparison
of the autocorrelation function of $U^{(h)}$ for the estimated and the parametric model gives the moving
average parameters of the MCARMA process.

\begin{proposition}
\label{propAutoCovUn}
Let $(U_n^{(h)})_{n\geq p}$ be the $d$-dimensional time series defined as $\Phi(B)Y_n^{(h)}=U_n^{(h)}$,
and $(\gamma_{U^{(h)}}(l))_{l\in\Z}=(\text{Cov}(U^{(h)}_{n+l},U^{(h)}_{n}))_{l\in\Z}$ denotes the autocovariance function.
%\begin{enumerate}
%\item[(a)]
Then  for $l=0,\ldots,p-1$:
\begin{eqnarray*}
\gamma_{U^{(h)}}(l)=\sum\limits_{\nu=1}^p\e^{hlR_\nu}\left[\sum\limits_{r=0}^{p-l-1}\sum_{\mu=1}^p \bigg(\e^{hrR_\nu}-\sum\limits_{j=1}^{r+l}\Phi_j\e^{h(r-j)R_\nu}\bigg)\Sigma_{\nu ,\mu}^{(h)}\bigg(\e^{h rR_\mu}-\sum\limits_{j=1}^{r}\Phi_j\e^{h(r-j)R_\mu}\bigg)^\mathsf{H}\right],
\end{eqnarray*}
and $\gamma_{U^{(h)}}(l)=0_{d\times d}$ for $l\geq p$, where
\begin{align*}
\Sigma_{\nu ,\mu}^{(h)}:=\Cov\left(N_{\nu,1}^{(h)} ,N_{\mu,1}^{(h)} \right)=\int_{0}^{h}\e^{R_\nu u}\res[F,R_\nu] ~\Sigma_L~\res[F,R_\mu]^\mathsf{H}{\e^{R_\mu^\mathsf{H}u}}\,\dif u.
\end{align*}
%\item[(b)] Define $\Theta_0:=I_d$ and $\Sigma_\varepsilon$ be the covariance matrix of the white noise $(\varepsilon_n^{(h)})_{n\in\Z}$
%given in \eqref{eqDiscrVARMA}. Then
%\begin{align*}
%\gamma_{U^{(h)}}(l)&=
%\left\{\begin{array}{ll}
%\sum_{k=0}^{p-|l|-1}\Theta_k\Sigma_\varepsilon\Theta_{k+|l|}^\mathsf{T},\quad &\text{ if }~~ |l|\leq p-1,
%\\
%0,  \quad &\text{ if }~~ |l|>p-1,
%\end{array}\right..
%\end{align*}
%\end{enumerate}
\end{proposition}
\begin{proof}
For lag $l\in\{0,\ldots,p-1\}$ we receive due to \eqref{eqDiscIIDSequence}:
\beao
\gamma_{U^{(h)}}(l)%\stackrel{\phantom{(3.33)}}{=}%&\Cov\left(U_{n+\vert l \vert}^{(h)},U_n^{(h)}\right)
%\notag\\
&\stackrel{\phantom{(3.33)}}{=}&\Cov\left(W_{0,n+ l}^{(h)}+\ldots+W_{p-1,n+l -p+1}^{(h)},W_{0,n}^{(h)}+\ldots+W_{p-1,n-p+1}^{(h)}\right)
\notag\\
&\stackrel{\phantom{(3.33)}}{=}&\sum_{r=0}^{p- l-1 }\Cov\left(W_{r+l,n-r}^{(h)},W_{r,n-r}^{(h)}\right)
\notag\\
&\stackrel{}{=}&\sum_{r=0}^{p-l-1 }\Bigg[\sum_{\nu=1}^p\sum_{\mu=1}^p
\bigg(\e^{h(r+ l )R_\nu}-\sum_{j=1}^{r+l }\Phi_j\e^{h(r+ l -j)R_\nu}\bigg)
\notag\\
&&\cdot \Cov\left(N_{\nu,n-r}^{(h)} ,N_{\mu,n-r}^{(h)} \right)  %\Cov\Bigg(\int\limits_{(n-i)h}^{(n-i+1)h}\e^{R_\nu[(n-i+1)h-u]}\res[F,R_\nu]\,\dif L(u),\int\limits_{(n-i)h}^{(n-i+1)h}\e^{R_\mu[(n-i+1)h-u]}\res[F,R_\mu]\,\dif L(u)\Bigg)
%\notag\\
\bigg(\e^{hrR_\mu}-\sum_{j=1}^{r}\Phi_j\e^{h(r-j)R_\mu} \bigg)^\mathsf{H}\Bigg],
\eeao
where
$
\Cov\left(N_{\nu,n-r}^{(h)} ,N_{\mu,n-r}^{(h)} \right)
=\Sigma_{\nu ,\mu}^{(h)},
$
and finally, the assertion follows.
\end{proof}

\end{document}